\documentstyle[11pt,dgspp]{article}
\newtheorem{theorem}{Theorem}[section]
\newtheorem{definition}[theorem]{Definition}
\newtheorem{proposition}[theorem]{Proposition}
\newtheorem{lemma}[theorem]{Lemma}

\newremark{example}[theorem]{Example}
\newremark{remark}[theorem]{Remark}

\newcommand{\A}{\Bbb A}
\newcommand{\B}{\Bbb B}

\newcommand{\bul}{\vcenter{\hbox{$\scriptscriptstyle\,\bullet$}}\,}

\newcommand{\ddu}{\ddot{u}}

\newcommand{\del}[1]{\nabla_{#1}}
\newcommand{\dg}{\dot\gamma}

\newcommand{\dps}{\displaystyle}
\newcommand{\ds}{\oplus} 

\newcommand{\du}{\dot{u}}

\newcommand{\dz}{\dot{z}}

\newcommand{\g}{\gamma}
\newcommand{\h}{\mbox{$\frak{h}$}}
\newcommand{\half}{\mbox{$\txs\frac{1}{2}$}}

\newcommand{\io}{\iota}
\newcommand{\Iaut}{I^{a\kern-.05em u\kern-.05em t}}
\newcommand{\Ispl}{I^{s\kern-.055em p\kern-.025em l}}
\newcommand{\jac}{R_{\,\dg}}

\newcommand{\lsp}{[\kern-0.15em[} 

\newcommand{\n}{\mbox{$\frak{n}$}}

\newcommand{\ph}{$p\kern-.1em H\!$}
\newcommand{\quar}{\mbox{$\txs\frac{1}{4}$}}
\newcommand{\rsp}{]\kern-0.15em]} 

\newcommand{\sixth}{\mbox{$\txs\frac{1}{6}$}}

\newcommand{\surj}{\rightarrow\kern-.82em\rightarrow}

\newcommand{\third}{{\txs\frac{1}{3}}}

\newcommand{\tquar}{\mbox{$\txs\frac{3}{4}$}}

\newcommand{\txs}{\textstyle}
\newcommand{\ve}{\mbox{$\varepsilon$}}

\newcommand{\z}{\frak{z}}

\newcommand{\II}{\rscr{I{\kern-.55em}I}}

\newcommand{\R}{\Bbb R}

\newcommand{\bZ}{\Bbb Z}

\renewcommand{\l}{\lambda}

\renewcommand{\v}{\frak v}

\makeatletter
\newcommand{\Ad}[1]{\mathop{\operator@font Ad}\nolimits_{#1}}
\newcommand{\ad}[1]{\mathop{\operator@font ad}\nolimits_{#1}}
\newcommand{\add}[2]{\mathop{\operator@font
   ad}\nolimits^{\dagger}_{#1}{\!#2}}
\newcommand{\Aut}{\mathop{\operator@font Aut}\nolimits}
\newcommand{\diag}{\mathop{\operator@font diag}\nolimits}
\newcommand{\End}{\mathop{\operator@font End}\nolimits}
\newcommand{\im}{\mathop{\operator@font im}\nolimits}
\newcommand{\mcp}{\mathop{\operator@font mult}\nolimits_{cp}}
\newcommand{\mev}{\mathop{\operator@font mult}\nolimits_{ev}}
\newcommand{\Ric}{\mathop{\operator@font Ric}\nolimits}
\newcommand{\rk}{\mathop{\operator@font rk}\nolimits}
\newcommand{\Dspecp}{\mathop{\operator@font {\escr D}spec}\nolimits_\wp}
\newcommand{\specl}{\mathop{\operator@font spec}\nolimits_\ell}
\newcommand{\specp}{\mathop{\operator@font spec}\nolimits_\wp}
\newcommand{\tr}{\mathop{\operator@font tr}\nolimits}
\makeatother

\preprint{JP4}
\title{Pseudo\boldmath$H$-type 2-step Nilpotent Lie Groups}
\author{Changrim Jang\thanks{On sabbatical leave; permanent address:
        Department of Mathematics,
        College of Natural Sciences,
        University of Ulsan,
        Ulsan 680-749,
        Republic of Korea.\protect\newline Email:
        crjang@mail.ulsan.ac.kr}
        \qquad\qquad$\!$Phillip E. Parker}
\address{Mathematics Department\\
        Wichita State University\\
        Wichita KS 67260-0033\\
        USA\\
        \hspace*{-1em}crjang@math.twsu.edu \qquad\qquad
        phil@math.twsu.edu}
\author{Keun Park}
\address{Department of Mathematics,\\
        College of Natural Sciences,\\
        University of Ulsan,\\
        Ulsan 680-749,\\
        Republic of Korea\\
        kpark@mail.ulsan.ac.kr}
\date{28 June 2003} 
\abstract{
Pseudo$H$-type is a natural generalization of $H$-type to geometries with
indefinite metric tensors.  We give a complete determination of the
conjugate locus including multiplicities.  We also obtain a partial 
characterization in terms of the abundance of totally geodesic,
3-dimensional submanifolds.
}
\msc{53C50}{22E25, 53B30, 53C30}

\begin{document}

\maketitle

\setcounter{page}{0}\thispagestyle{empty}\strut\vfill\eject

\section{Introduction}
\label{intro}

The study of the conjugate locus in Riemannian Lie groups of $H$-type has 
evolved in parallel with that in the ambient class of Riemannian 2-step 
nilpotent Lie groups. For a brief history of the latter see the 
introduction to \cite{JP1}; here we shall give some for the groups of
$H$-type.

Groups of $H$-type were introduced by Kaplan in 1981 \cite{K,K1}.  A
little later, Boggino in 1985 \cite{B} found many of the conjugate points
in a Reimannian group of $H$-type.  In 1995, Berndt, Tricerri and Vanhecke
\cite{BTV} determined almost all of the conjugate points in Riemannian
groups of $H$-type, and all of them when the center is 1-dimensional.
They included determination of the multiplicities in most cases.  Walschap
in 1997 \cite{W} showed that for Riemannian $H$-type groups, the cut locus
is the conjugate locus.  Finally, in unpublished work from 2001, Kim and
Park \cite{KP} found explicit formulas for all the conjugate points in any
Riemannian $H$-type group.  A recent paper \cite{JK} determined all
conjugate points in quaternionic Heisenberg groups with the Killing metric
tensor from $Sp(1,n)$.

Using Ciatti's generalization of $H$-type to pseudo$H$-type \cite{C}, in 
this paper we shall determine explicit formulas for all conjugate points 
and their multiplicities in the (larger) class of pseudo$H$-type, 2-step
nilpotent Lie groups with a left-invariant pseudoriemannian (indefinite)
metric tensor.  When specialized to positive-definite metric tensors, we
recover all the known Riemannian results for $H$-type groups with shorter,
simpler, more conceptual proofs.

\medskip
By an {\em inner product\/} on a vector space $V$ we shall mean a
nondegenerate, symmetric bilinear form on $V$, generally denoted by
$\langle\,,\rangle$.  Our convention is that $v\in V$ is timelike if
$\langle v,v\rangle > 0$, null if $\langle v,v\rangle = 0$, and spacelike
if $\langle v,v\rangle < 0$.

Throughout, $N$ will denote a connected, 2-step nilpotent Lie group with
Lie algebra \n\ having center $\z$.  We shall use $\langle\,,\rangle$ to
denote either an inner product on \n\ or the induced left-invariant
pseudoriemannian (indefinite) metric tensor on $N$. We also assume that
the center $\z$ is nondegenerate with respect to such a metric tensor.

We denote the adjoint with respect to $\langle\,,\rangle$ of the adjoint
representation of the Lie algebra \n\ on itself by $\add{}{}\!$.  In the
case of a nondegenerate center, the involution $\iota$ of \cite{CP4} is
merely given by $\io (z_\alpha )=\ve_\alpha\, z_\alpha$ and $\io (e_a )=
\bar{\ve}_a \, e_a$ where $\langle z_\alpha ,z_\alpha \rangle =
\ve_\alpha$ and $\langle e_a,e_a \rangle = \bar{\ve}_a$ on an orthonormal
basis of \n.  Then the operator $j:\z\rightarrow\End\left(\v\right) $ is
given by $ j(z)x = \io\add{x}{\io z}$.  We refer to \cite{CP4} and
\cite{O} for other notations and results.

For convenience, we shall use the notation $J_z = \add{\bul}{z}$ for any
$z\in\z$. (Since the center is nondegenerate, the involution $\io$ may
be omitted.) We follow \cite{C} for this next definition. As in the 
Riemannian case, one might as well make 2-step nilpotent part of the 
definition since it effectively is so anyway.
\begin{definition}
$N$ is\label{pht} said to be of\/ {\em pseudo$H$-type} if and only if
$$ J^2_z = -\langle z, z \rangle I $$
for any $z \in \z$.
\end{definition}
Lemma \ref{l3} collects some basic properties of these spaces.
In Theorem \ref{th1} we complete the determination of the conjugate locus
for such groups; compare with Theorem 2.2 in \cite{JP1}.  Note that, in
contrast to Theorem 2.4 and Corollary 2.6 in \cite{JP1}, the center may
now be of dimension greater than 1.

For $H$-type groups, the {\em nonsingular\/} condition was shown to be an
important determining property of their geometries by Eberlein \cite{E}.
For pseudo$H$-type groups, however, it must be replaced by something
better adapted to pseudoriemannian (indefinite metric tensor) geometries.
Indeed, one may easily show that a nonsingular, pseudo$H$-type group
necessarily has a (positive- or negative-) definite center.
\begin{definition}
Let\label{pns} $N$ be a pseudoriemannian, simply connected, 2-step nilpotent
Lie group with Lie algebra $\n$ with nondegenerate center $\z$.  We say
that $N$ is\/ {\em pseudoregular} if and only if (1) $\ad{x}$ is
surjective for every nonnull $x\in\v=\z^\perp$, and (2) $J_z$ is
nonsingular for every nonnull $z\in\z$.
\end{definition}
Although (1) and (2) are equivalent in the Riemannian (definite) case, it
is easy to show by examples that they are independent in general. It is
not clear how this definition should be extended to handle a degenerate
center. Similarly to the Riemannian case in which $H$-type implies 
nonsingular, it is easy to see that pseudo$H$-type implies pseudoregular.

In Theorem \ref{th2} we generalize Eberlein's Theorem 6.1 from \cite{E} to
pseudoregular groups of pseudo$H$-type, obtaining essentially the same
partial characterization of pseudo$H$-type in terms of the abundance of
totally geodesic submanifolds as Eberlein did there of $H$-type.

The rest of this Section consists of preliminaries and recollections. In
Section \ref{mr} we state the main results and state and prove some
secondary results. Sections \ref{pf1} and \ref{pf2} are devoted to the
proofs of the main Theorems \ref{th1} and \ref{th2}, respectively.

\medskip
We begin by specializing Theorems 3.1 and 3.6 of \cite{CP4} to
the case of a nondegenerate center.
\begin{proposition}
For\label{conn} all $z,z'\in\z$ and $e,e'\in\v$ we have
\begin{eqnarray*}
\del{z}z' & = & 0\, ,\\
\del{z}e = \del{e}z &=& -\half J_z e\, ,\\
\del{e}e' &=& \half [e,e']\,.
\end{eqnarray*}
\end{proposition}
\begin{proposition}
For\label{curv} all $z,z',z''\in\z$ and $e,e',e''\in\v$ we have
\begin{eqnarray*}
R(z,z')z'' &=& 0\, ,\\
R(z,z')e   &=& \quar (J_z J_{z'} e - J_{z'} J_z e)\, ,\\
R(z,e)z'   &=& \quar J_z J_{z'} e\, ,\\
R(z,e)e'   &=& \quar [e, J_z e']\, ,\\
R(e,e')z   &=& -\quar\big( [e, J_z e'] +[ J_z e,e']\big)\, ,\\
R(e,e')e'' &=& \quar\big( J_{[e,e'']}e' - J_{[e',e'']}e\big)
               +\half J_{[e,e']}e''\,.
\end{eqnarray*}
\end{proposition}

To study conjugate points, we use the Jacobi operator.
\begin{definition}
Along\label{jo} the geodesic $\g$, the\/ {\em Jacobi operator} is given by
$$ \jac \bul = R(\bul,\dg)\dg\,. $$
\end{definition}
In physics, this operator measures the relative acceleration produced by
tidal forces along $\g$ \cite[p.\,219]{O}.  For the reader's convenience,
we recall that a {\em Jacobi field\/} along $\g$ is a vector field along
$\g$ which is a solution of the {\em Jacobi equation}
$$ \nabla_{\dg}^2 Y(t)+\jac Y(t) = 0 $$
along $\g$. The point $\g(t_0)$ is {\em conjugate\/} to the point $\g(0)$ 
if and only if there exists a nontrivial Jacobi field $Y$ along $\g$ such 
that $Y(0) = Y(t_0) = 0$.

Next, we specialize Proposition 4.8 of \cite{CP4} to our present setting.
As there, $L_n$ denotes left translation in $N$ by $n\in N$.
\begin{proposition}
Let\label{e3.2} $N$ be simply connected, $\g$ a geodesic with $\g(0)=1\in N$,
and $\dot{\g}(0) = z_0 + x_0 \in \n$. Then
$$\dot{\g}(t) = L_{\gamma(t)*}\left( z_0 + e^{tJ}x_0 \right) $$
where $J = \add{\bul}{z_0}$.
\end{proposition}
As in \cite{JPk}, we shall identify all tangent spaces with $\n = T_1 N$.
Thus we regard
$$ \dg(t) = z_0 + e^{tJ}x_0 $$
as the geodesic equation. Using Proposition \ref{curv}, we compute
\begin{lemma}
The\label{jof} Jacobi operator along the geodesic $\g$ in $N$ with $\g(0)
= 1$ and $\dg(0) = z_0 + x_0$ is given by
\begin{eqnarray*}
\jac(z+x) &=& \tquar J_{[x,x']}x' + \half J_z Jx' - \quar
   JJ_z x' - \quar J^2 x \\[.5ex]
&& {} - \half[x, Jx'] + \quar[x',Jx] + \quar[x',J_z x']
\end{eqnarray*}
for all $z\in\z$ and $x\in\v$, where $x' = e^{tJ}x_0$ and $J=J_{z_0}$.\eop
\end{lemma}

For the reader's convenience, we provide the statement of Proposition 2.1
from \cite{JP1}.
\begin{proposition}
Let\label{yj} $\g$ be a geodesic with $\g(0) =1$ and $\dg(0)
= z_0 + x_0 \in \z\ds\v=\n$.  A vector field $Y(t) = z(t) + e^{tJ}u(t)$
along $\g$, where $ z(t) \in \z $ and $u(t) \in \v$ for each $t$, is a
Jacobi field if and only if
\begin{eqnarray*}
\dot{z}(t) -[e^{tJ}u(t) , x'(t)] &=& \zeta\, , \\
e^{tJ}\ddot{u}(t)+e^{tJ}J\dot{u}(t)-\add{x'(t)}{\zeta } &=& 0\, ,
\end{eqnarray*}
where $x'(t) = e^{tJ}x_0$ with $J=J_{z_0}$ and $\zeta\in\z$ is a constant.
\end{proposition}
We also provide an adapted version of the statement of Theorem 2.2 from
that paper. First we recall some notation.
\begin{definition}
Let $\g$ denote a geodesic and assume that $\g(t_0)$ is conjugate to 
$\g(0)$ along $\g$. To indicate that the multiplicity of $\g(t_0)$ is $m$,
we shall write $\mcp(t_0)=m$. To distinguish the notions clearly, we shall
denote the multiplicity of $\l$ as an eigenvalue of a specified linear 
transformation by $\mev{\l}$.
\end{definition}
Let $\g$ be a geodesic with $\g(0)=1$ and $\dg(0) = z_0+x_0 \in
\z\ds\v$, respectively, and let $J = J_{z_0}$. If $\g$ is not
null, we may assume $\g$ is normalized so that $\langle\dg,\dg\rangle =
\pm 1$. As usual, $\bZ^*$ denotes the set of all integers with 0 removed.
\begin{theorem}
Under\label{th2.2} these assumptions, if $N$ is of pseudo$H$-type then:
\begin{enumerate}
\item if $z_0=0$ and $x_0\neq 0$, then $\g(t)$ is conjugate to $\g(0)$
along $\g$ if and only if $\langle x_0,x_0\rangle < 0$ and
$$-\frac{12}{t^2} = \langle x_0,x_0\rangle\, ,$$
in which case $\mcp(t) = \dim\z$;
\item if $z_0\neq 0$ and $x_0 = 0$, then $\g(t)$ is conjugate to $\g(0)$
along $\g$ if and only if $\langle z_0,z_0\rangle > 0$ and
$$t \in  \frac{2\pi}{|z_0|} \bZ^* ,$$
in which case $\mcp(t) = \dim\v$.
\end{enumerate}
\end{theorem}
Note that Parts 1 and 2 here correspond with Parts 2 and 3 there. Since we
are assuming $N$ is of pseudo$H$-type, $J_z=0$ only for $z=0$. Indeed,
writing $z = z_t-z_s$ where $z_t$ is timelike and $z_s$ is spacelike, it
is immediate that $J_{z_t}^2\neq J_{z_s}^2$ whence $J_{z_t}\neq J_{z_s}$
so $J_z\neq 0$ unless $z=0$.  Thus Part 1 of Theorem 2.2 in \cite{JP1}
cannot occur here, and we have restated the other two parts appropriately.

Finally, here are some basic results on pseudo$H$-type spaces, useful for
computations in them.
\begin{lemma}
Let\label{l3} $N$ be a group of pseudo$H$-type. Then these hold for all
$z, z' \in \z$ and $x, y\in\v$:
\begin{enumerate}
\item $\langle J_zx,J_{z'}x\rangle = \langle z,z'\rangle \langle
   x,x\rangle $;
\item $\langle J_zx,J_zy\rangle = \langle z,z\rangle\langle x,y\rangle$;
\item $J_zJ_{z'}+J_{z'}J_z = -2\langle z,z'\rangle I $;\label{i3}
\item $[x,J_zx] = \langle x,x\rangle z$\,.\eop
\end{enumerate}
\end{lemma}

\section{Main Results}
\label{mr}

First, we present an example of a pseudo$H$-type group which is singular
(not nonsingular), therefore not of $H$-type.
\begin{example}
Let $\n$ be the following Lie algebra with indefinite form
$\langle\,,\rangle$. We present it {\em via\/} a null basis for $\v$ and
an orthonormal basis for $\z$.
$$ \n = \lsp x_i,y_i,v_i,w_i \mid 1\le i\le k \rsp\ds\lsp z_1,z_2,z_3 \rsp
$$
where all the nontrivial inner products are given by
$$ \langle x_i,w_i \rangle = -1\,, \qquad \langle y_i,v_i\rangle = 1 $$
$$ \langle z_1,z_1\rangle = 1 = -\langle z_2,z_2\rangle = -\langle
z_3,z_3\rangle\,. $$
The structure equations are
$$ [x_i,v_i] = \half(z_1-z_2)\, ,\qquad [x_i,w_i] = z_3\, , $$
$$ [y_i,v_i] = z_3\, ,\qquad [y_i,w_i] = 2(z_1+z_2)\, , $$
with all other brackets vanishing.
The group $N$ is the unique, simply connected, 2-step nilpotent Lie group
with Lie algebra $\n$.  It is of pseudo$H$-type:  a straightforward
computation shows that $J_{az_1+bz_2+cz_3}^2 = -(a^2-b^2-c^2)I$ as
required. It is easy to verify that it is singular.
\end{example}

The main difference between pseudo$H$-type and $H$-type is that $J_z$ can
drop rank on null vectors. This can be quantified.
\begin{proposition}
Let $N$ be of pseudo$H$-type with Lie algebra $\n=\z\ds\v$. Then $\dim\v$
is even and $\rk J_z = \half\dim\v$ for nonzero null $z$.
\end{proposition}
\begin{proof}
Choose a nonnull $z\in\z$. First, assume $\langle z,z\rangle=\l^2>0$, so
$J_z^2=-\l^2I$. Now choose a unit (timelike or spacelike) $e_1\in\v$; then
$e_1\perp J_ze_1$. Next, choose a unit (timelike or spacelike) $e_2\in\v$
such that $\langle e_1,e_2\rangle=\langle e_2,J_ze_2\rangle=0$; then
$e_1$, $J_ze_1$, $e_2$, $J_ze_2$ are mutually orthogonal. Continuing this
process, we obtain a basis for $\v$ with an even number of vectors.

Now assume $\langle z,z\rangle=-\l^2<0$, so $J_z^2=\l^2I$. This implies
that $\v=\ker(J_z+\l I)\ds\ker(J_z-\l I)$. One may show that $\ker(J_z+\l
I)$ and $\ker(J_z-\l I)$ are complementary null subspaces; {\em cf.}\
\cite{JP2} for such an argument. Again, we find $\dim\v$ is even.

For the rank of $J_z$, let $0\neq z\in\z$ be null and assume that
$J_zx=0$.  Write $z=z_t-z_s$ where $z_t$ is timelike, $z_s$ is spacelike,
and $z_t\perp z_s$.  By assumption, $J_{z_t}x=J_{z_s}x$.  By Lemma
\ref{l3}, $\langle J_{z_t}x,J_{z_s}x\rangle = \langle z_t,z_s\rangle
\langle x,x\rangle = 0$ on the one hand, and $\langle
J_{z_t}x,J_{z_s}x\rangle = \langle J_{z_t}x,J_{z_t}x\rangle = \langle
z_t,z_t\rangle \langle x,x\rangle$ on the other.  Hence $x$ is null and
$\ker J_z$ is a null subspace of $\v$.  Thus there exists a complementary
null subspace $\v'$.  Consider the subspace $\ker J_z\ds\v'$.  Note that
$J_z\v\subseteq\ker J_z$ since $J_z^2=0$.  Now, $J_z(\ker J_z\ds\v') =
J_z\v'$ and $\dim J_z\v' = \dim\v' = \dim\ker J_z$ so $J_z\v' = J_z\v =
\ker J_z$.  Therefore $\v=\ker J_z\ds\v'$ and $\rk J_z = \dim\v' =
\half\dim\v$.
\end{proof}

Let $\g$ be a geodesic with $\g(0)=1$ and $\dg(0) = z_0+x_0 \in \z\ds\v$,
respectively, and let $J = J_{z_0} = \add{\bul}{z_0}$.  If $\g$ is not
null, we may assume $\g$ is normalized so that $\langle\dg,\dg\rangle =
\pm 1$.  As usual, $\bZ^*$ denotes the set of all integers with 0 removed.
\begin{theorem}
Let\label{th1} $\g$ be such a geodesic in a pseudo$H$-type group $N$ with
$z_0 \neq 0 \neq x_0$.
\begin{enumerate}
\item If $\langle z_0 , z_0 \rangle =\alpha^2$ with $\alpha>0$, then
$\gamma (t_0)$ is conjugate to $\g(0)$ along $\gamma$ if and only if
$$ t_0 \in \frac{2\pi}{\alpha }\bZ^* \cup \A_1 \cup \A_2$$
where
$$ \A_1 = \left\{ t \in \R \biggm| \langle x_0 , x_0 \rangle\frac{\alpha
   t}{2}\cot\frac{\alpha t}{2} = \langle \dg , \dg \rangle \right\}
   \quad\mbox{and \qquad } $$
$$ \A_2 = \left\{ t\in \R \biggm| \alpha t = \frac{\langle x_0 , x_0
   \rangle }{\langle \dg , \dg \rangle + \langle z_0 , z_0 \rangle }
   \sin\alpha t\right\}\quad\mbox{when }\,\dim\z  \geq 2\,. $$
If $t_0\in (2\pi/\alpha)\bZ^*$, then
$$ \mcp(t_0) = \left\{ \begin{array}{cl}
\dim\v-1 & \mbox{ if }\,\langle\dg,\dg\rangle + \langle z_0,z_0\rangle \ne
   0\,,\\
\dim\n-2 & \mbox{ if }\,\langle\dg,\dg\rangle + \langle z_0,z_0\rangle = 
   0\,. \end{array} \right. $$
If $t_0\notin (2\pi/\alpha)\bZ^*$, then
$$ \mcp(t_0) = \left\{ \begin{array}{cl}
1 & \mbox{ if }\, t_0\in\A_1 - \A_2\,,\\
\dim\z-1 & \mbox{ if }\, t_0\in\A_2 - \A_1\,,\\
\dim\z & \mbox{ if }\, t_0\in\A_1\cap\A_2\,. \end{array} \right. $$

\item If $\langle z_0 , z_0 \rangle = -\beta^2$ with $\beta >0$ , then
$\gamma (t_0)$ is a conjugate point along $\gamma$ if and only if $t_0 \in
\B_1 \cup \B_2$ where
$$ \B_1 = \left\{ t \in \R \biggm| \langle x_0 , x_0 \rangle \frac{\beta
   t}{2}\coth\frac{\beta t}{2} = \langle \dg , \dg \rangle \right\}
   \quad\mbox{and \qquad } $$
$$ \B_2 = \left\{ t\in \R \biggm| \beta t = \frac{\langle x_0 , x_0
   \rangle }{\langle \dg , \dg \rangle + \langle z_0 , z_0 \rangle } \sinh
   \beta t\right\}\quad\mbox{when }\,\dim\z \geq 2\,. $$
The multiplicity is
$$ \mcp(t_0) = \left\{ \begin{array}{cl}
1 & \mbox{ if }\, t_0\in\B_1 - \B_2\,,\\
\dim\z-1 & \mbox{ if }\, t_0\in\B_2 - \B_1\,,\\
\dim\z & \mbox{ if }\, t_0\in\B_1\cap\B_2\,. \end{array} \right. $$

\item If $ \langle z_0 , z_0 \rangle = 0$, then $ \gamma (t_0) $ is a
conjugate point along $\gamma$ if and only if
$$ t_0^2 = -\frac{12}{\langle x_0 , x_0 \rangle }\, ,$$
and $\mcp(t_0)=\dim\z-1$.

\end{enumerate}
\end{theorem}
Together with Theorem \ref{th2.2} \cite[2.2]{JP1}, this covers all
cases for a pseudo$H$-type group with a center of any dimension. The proof
occupies Section \ref{pf1}.

Here is our generalization to pseudoregular groups of Eberlein's theorem
for nonsingular groups \cite[Thm.\,6.1]{E}.  We note that O'Neill
\cite[Ex.\,9, p.\,125]{O} has extended the definition of totally geodesic
to (possibly) degenerate submanifolds of pseudoriemannian manifolds.  We
shall use this extended version.
\begin{theorem}
Let\label{th2} $N$ be pseudoregular.  For every geodesic $\g$ with
$\g(0)=1\in N$ and $\dg(0)=z_0+x_0$ with $|z_0|\neq 0$ there exists a
connected, totally geodesic, 3-dimensional submanifold $H$ such that
$\dg(0)\in \h=T_1H$ and $\h\cap\z$ is nonnull, if and only if the metric
tensor on $N$ is homothetic to one of pseudo$H$-type.
\end{theorem}
The proof is found in Section \ref{pf2}.

To show that pseudoregular is necessary, here is a group that satisfies
condition (2), but not condition (1), of Definition \ref{pns} and
the hypothesis {\em supra\/} about totally geodesic submanifolds, but is
not homothetic to one of pseudo$H$-type.
\begin{example}
Let $N$ be the unique, simply connected, 2-step nilpotent Lie group with
Lie algebra $\n = \lsp e_1,e_2\rsp\ds\lsp z_1,z_2\rsp$, structure equation
$[e_1,e_2]=z_1-z_2$, and nontrivial inner products
$$ \langle e_1,e_1\rangle = \langle e_2,e_2\rangle = 1\, ,$$
$$ \langle z_1,z_1\rangle = -\langle z_2,z_2\rangle = 1\,. $$
One readily computes $J^2_{z_1+z_2}=-4I$, but $z_1+z_2$ is null so $N$ is
not of pseudo$H$-type, nor even homothetic to one.

It is easy to see that it satisfies condition (2) of Definition \ref{pns},
and $\ad{e_1}$ is clearly not surjective so condition (1) fails to hold.
Thus $N$ is not pseudoregular.
With a bit more effort, one may show that for every geodesic $\g$ with
$\g(0)=1\in N$ there exists a 3-dimensional, totally geodesic submanifold
$H$ such that $\dg(0)\in T_1H$.  Indeed, if $\dg(0)=z_0+x_0$ then $H =
\{\exp_1(ax_0+bz_1+cz_2)\mid a,b,c\in\R\}$ where we use the
pseudoriemannian (geometric) exponential map, {\em not\/} the Lie group
one.  In fact, such an $H$ is totally geodesic for any $x_0\in\v$.
\end{example}

\section{Proof of Theorem \protect\ref{th1}}
\label{pf1}

For the first part, assume $ Y(t) = z(t) + e^{tJ}u(t)$ is a nontrivial
Jacobi field along $\gamma$ with $ Y(0)=Y(t_0 )= 0$. Then by Proposition 
\ref{yj}, we may assume that
\begin{eqnarray}
\dz(t) - [ e^{tJ} u(t) , x'(t)] = cz_0 + \zeta\label{ej1}\\
e^{tJ} \ddu(t) + e^{tJ}J \du(t) - J_{cz_0 +\zeta}x'(t) = 0\label{ej2}
\end{eqnarray}
for a constant $c$ and a constant vector $ \zeta \in\z$ with $ \langle
z_0, \zeta \rangle = 0$.

The general solution of equation (\ref{ej2}) satisfying $u(0)=0$ is
\begin{equation}
\label{ej3}
u(t) = ctx_0 + (e^{-tJ} -I )v_0 - \frac{1}{2 \langle z_0 , z_0
   \rangle }(e^{-2tJ}-e^{-tJ} )J_\zeta x_0
\end{equation}
for some $v_0\in\v$.  To get this, we used the fact that $ e^{-tJ}J_\zeta
= J_\zeta e^{tJ}$; this follows from item \ref{i3} in Lemma \ref{l3}.

Assume that
\begin{equation}
\label{hyp1}
t_0 \in \frac{2\pi}{|z_0|}\bZ^*.
\end{equation}
Then from (\ref{ej3}) we have $u(t_0) = ct_0x_0 = 0$ whence $c=0$. Thus
\begin{equation}
\label{new4}
u(t) = (e^{-tJ} -I )v_0 - \frac{1}{2 \langle z_0 , z_0
   \rangle }(e^{-2tJ}-e^{-tJ} )J_\zeta x_0\, ,
\end{equation}
and from (\ref{ej1}) and (\ref{new4})
$$ \dz(t) + \big[ (e^{tJ}-I)v_0 , x'(t) \big] -\frac{1}{2\langle z_0 , z_0
   \rangle }\big[ (e^{-tJ}-I)J_\zeta x_0 , x'(t)\big] = \zeta\,.$$
Using
$$ e^{tJ} = \left(\cos|z_0|t\right)I + \frac{1}{|z_0|}\left(
   \sin|z_0|t\right)J\, , $$
this becomes\arraycolsep .1em
\begin{eqnarray*}
\dz(t) &-& \left[ v_0, \cos|z_0|t\,x_0 + \frac{1}{|z_0|}\sin|z_0|t\,Jx_0 
   \right] \\
&+&\left[ \cos|z_0|t\,v_0 + \frac{1}{|z_0|}\sin|z_0|t\,Jv_0, 
   \cos|z_0|t\,x_0 + \frac{1}{|z_0|}\sin|z_0|t\,Jx_0 \right] \\
&+&\frac{1}{2\langle z_0,z_0\rangle }\,\bigg[ \cos|z_0|t\,J_\zeta x_0 - 
   \frac{1}{|z_0|}\sin|z_0|t\,JJ_\zeta x_0, \\
&& \hspace{10em}\cos|z_0|t\,x_0 + \frac{1}{|z_0|}\sin|z_0|t\,Jx_0 \bigg] 
   \\
&-&\frac{1}{2\langle z_0,z_0\rangle } \left[ J_\zeta x_0, \cos|z_0|t\,x_0 
   + \frac{1}{|z_0|}\sin|z_0|t\,Jx_0 \right] = \,\zeta\,.
\end{eqnarray*}
Also from Lemma \ref{l3}, we get $[J_\zeta x_0,Jx_0]=[JJ_\zeta x_0,x_0]$
and $[JJ_\zeta x_0,Jx_0] = \langle x_0,x_0\rangle\langle
z_0,z_0\rangle\zeta$.  Together with some computation, these yield
\begin{eqnarray*}
\dz(t) &+& \left(\cos^2|z_0|t - \cos|z_0|t\right)[v_0,x_0] \\
&+&\frac{1}{|z_0|}\left(\cos|z_0|t\sin|z_0|t - \sin|z_0|t\right) 
   [v_0,Jx_0] \\
&+& \frac{1}{|z_0|} \sin|z_0|t\cos|z_0|t\,[Jv_0,x_0] 
   + \frac{1}{|z_0|^2} \sin^2|z_0|t\,[Jv_0,Jx_0] \\
&+& \frac{1}{2\langle z_0,z_0\rangle } \left(\cos|z_0|t - 
   \cos^2|z_0|t\right) \langle x_0,x_0\rangle\zeta \\
&-& \frac{1}{2\langle z_0,z_0\rangle } \sin|z_0|t\,[J_\zeta x_0,Jx_0] 
   - \frac{1}{2\langle z_0,z_0\rangle } \sin^2|z_0|t\langle 
   x_0,x_0\rangle\zeta \,=\, \zeta\,.
\end{eqnarray*}
Integrating under $z(0)=0$,
\begin{eqnarray}
z(t) &+& \left( \frac{t}{2} +\frac{1}{4|z_0|}\sin 2|z_0|t - 
   \frac{1}{|z_0|}\sin|z_0|t\right)[v_0,x_0] \nonumber\\
&+& \frac{1}{|z_0|^2} \left(\cos|z_0|t - \quar\cos 2|z_0|t - 
   \tquar\right)[v_0,Jx_0] \nonumber\\
&+& \frac{1}{|z_0|^2} \left(\frac{t}{2} - \frac{1}{4|z_0|}
   \sin2|z_0|t\right)[Jv_0,Jx_0] \nonumber\\
&+& \frac{1}{4|z_0|^2} \left(1-\cos^2|z_0|t\right)[Jv_0,x_0]
   + \frac{\langle x_0,x_0\rangle }{2\langle z_0,z_0\rangle } 
   \left(\frac{1}{|z_0|}\sin|z_0|t - t\right)\zeta \nonumber\\
&+& \frac{1}{2|z_0|^3}\left(\cos|z_0|t-1\right)[J_\zeta x_0,Jx_0] \,=\, 
   t\zeta\,. \label{new6}
\end{eqnarray}\arraycolsep .2em
If $\langle z_0,z'\rangle = 0$ for some $z'\in\z$, then
\begin{eqnarray*}
\langle[v_0,x_0] + \frac{1}{|z_0|^2}[Jv_0,Jx_0],z'\rangle
&=& \langle[v_0,x_0],z'\rangle + \frac{1}{|z_0|^2} \langle[Jv_0,
   Jx_0],z'\rangle \\
&=& \langle J_{z'}v_0,x_0\rangle + \frac{1}{|z_0|^2} \langle[J_{z'}Jv_0, 
   Jx_0\rangle \\
&=& \langle J_{x'}v_0,x_0\rangle - \frac{1}{|z_0|^2} \langle[JJ_{z'}v_0, 
   Jx_0\rangle \\
&=& \langle J_{z'}v_0,x_0\rangle - \langle J_{z'}v_0,x_0\rangle \\
&=& 0\,.
\end{eqnarray*}
This implies that
$$ [v_0,x_0] + \frac{1}{|z_0|^2}[Jv_0,Jx_0] = -2\frac{\langle 
   v_0,Jx_0\rangle }{\langle z_0,z_0\rangle }\, z_0\,. $$
Substituting into (\ref{new6}),
\begin{eqnarray}
\lefteqn{\hspace{-1.5em}z(t) \,=\, } \nonumber\\
&&  \left( \frac{1}{|z_0|}\sin|z_0|t - \frac{1}{4|z_0|}\sin 2|z_0|t 
   \right)[v_0,x_0] \nonumber\\
&&{}+\frac{1}{|z_0|^2} \left( \quar\cos 2|z_0|t - \cos|z_0|t +
   \tquar\right)[v_0,Jx_0] \nonumber\\
&&{}+\frac{1}{4|z_0|^2} \left(\cos^2|z_0|t-1\right)[Jv_0,x_0] 
   + \frac{1}{4|z_0|^3} \sin2|z_0|t[Jv_0,Jx_0] \nonumber\\
&&{}+\left( \frac{\langle x_0,x_0\rangle+2\langle z_0,z_0\rangle 
   }{2\langle z_0,z_0\rangle }\,t - \frac{\langle x_0,x_0\rangle 
   }{2|z_0|^3}\sin|z_0|t\right)\zeta 
   +\frac{\langle v_0,Jx_0\rangle }{\langle z_0,z_0\rangle }t\,z_0 
   \quad\label{new7}
\end{eqnarray}
Since $z(t_0)=0$, this gives
$$ t_0\left(\frac{\langle v_0,Jx_0\rangle }{\langle z_0,z_0\rangle }\,z_0 
+ \frac{\langle x_0,x_0\rangle+2\langle z_0,z_0\rangle }{2\langle
z_0,z_0\rangle }\,\zeta\right) = 0 $$
which implies that $\langle v_0,Jx_0\rangle = 0$ and $\left( \langle 
x_0,x_0\rangle+2\langle z_0,z_0\rangle \right)\zeta=0$.

\smallskip
If $\langle x_0,x_0\rangle+2\langle z_0,z_0\rangle \neq 0$, then 
$\zeta=0$. Now (\ref{new4}) and (\ref{ej1}) yield $u(t) = \left( 
e^{tJ}-I\right)v_0$ with $\langle v_0,Jx_0\rangle=0$, so (\ref{new7}) 
simplifies to
\begin{eqnarray*}
z(t) &=& \left( \frac{1}{|z_0|}\sin|z_0|t - \frac{1}{4|z_0|}\sin 2|z_0|t 
   \right)[v_0,x_0] \\
&&{}+\frac{1}{|z_0|^2} \left( \quar\cos 2|z_0|t - \cos|z_0|t +
   \tquar\right)[v_0,Jx_0] \\
&&{}+\frac{1}{4|z_0|^2} \left(\cos^2|z_0|t-1\right)[Jv_0,x_0] 
   + \frac{1}{4|z_0|^3} \sin2|z_0|t[Jv_0,Jx_0]\,.
\end{eqnarray*}
It readily follows that $\mcp(t_0) = \dim\v-1$ as claimed.

\smallskip
If $\langle x_0,x_0\rangle+2\langle z_0,z_0\rangle = 0$, then similarly
$$ u(t) = \left(e^{tJ}-I\right)v_0 - \frac{1}{2\langle z_0,z_0\rangle } 
   \left(e^{-2tJ}-e^{-tJ}\right)J_\zeta x_0\, , $$
\begin{eqnarray*}
z(t) &=& \left( \frac{1}{|z_0|}\sin|z_0|t - \frac{1}{4|z_0|}\sin 2|z_0|t 
   \right)[v_0,x_0] \\
&&{}+\frac{1}{|z_0|^2} \left( \quar\cos 2|z_0|t - \cos|z_0|t +
   \tquar\right)[v_0,Jx_0] \\
&&{}+\frac{1}{4|z_0|^2} \left(\cos^2|z_0|t-1\right)[Jv_0,x_0] 
   + \frac{1}{4|z_0|^3} \sin2|z_0|t[Jv_0,Jx_0] \\
&&{}- \frac{\langle x_0,x_0\rangle }{2|z_0|^3}\sin|z_0|t\,\zeta\, ,
\end{eqnarray*}
and it follows that $\mcp(t_0) = \dim\v-1 + \dim\z-1 = \dim\n-2$ as 
desired.

\medskip
Now assume that
\begin{equation}
\label{new8}
t_0 \notin \frac{2\pi}{|z_0|}\bZ^*.
\end{equation}
If $ cz_0 + \zeta = 0$ in (\ref{ej2}), then we have $ u(t) =
(e^{-tJ}-I)v_0$ for a constant vector $v_0 \in \v$.  Since $u(t_ 0) =
(e^{-t_0 J}-I)v_0 = 0$, by our assumption (\ref{new8}) we must have
$v_0=0$, which implies that $u(t)=0$.  This and (\ref{ej1}) with $z(0)=0$
imply that $z(t) =0$ for all $t$.  Thus $Y= 0$, contradicting the
nontriviality of $Y$.  Therefore $cz_0 + \zeta \neq 0$.

{}From (\ref{ej3}), $u(t_0)=0$, and (\ref{new8}), 
\begin{equation}
\label{new9}
\begin{array}{rcl}
\dps u(t) &=&\dps ctx_0 - (e^{-tJ} -I )(e^{-t_0 J}-I)^{-1}ct_0 x_0
   \\[1ex]
&&\dps {}+\frac{1}{2 \langle z_0 , z_0 \rangle }(e^{-tJ}-I )
   (e^{-t_0 J} -e^{-tJ} )J_\zeta x_0\, .\end{array}
\end{equation}
Substituting into equation (\ref{ej1}) we get
\begin{eqnarray*}
\lefteqn{\dz(t) + \big[ (I- e^{tJ})(e^{-t_0 J}- I)^{-1}ct_0 x_0 , x'(t)
   \big]}\\[.5ex]
&& {}-\frac{1}{2\langle z_0 , z_0 \rangle }\big[ (I-e^{tJ})(e^{-t_0
   J}-e^{-tJ})J_\zeta x_0 , x'(t)\big]\, =\, cz_0 +\zeta\,.
\end{eqnarray*}
Using the identities
\begin{eqnarray*}
e^{tJ} &=& \left(\cos|z_0|t\right)I + \frac{1}{|z_0|}\left(
   \sin|z_0|t\right)J\, ,\\
\left(e^{-t_0J}-I\right)^{-1} &=& -\half I +
   \frac{1}{2|z_0|}\cot\frac{|z_0|t_0}{2} J\, ,
\end{eqnarray*}
this becomes
\begin{eqnarray*}
\dz(t) + \frac{ct_0}{2|z_0|} \left( \cot\frac{|z_0|t_0}{2} - \sin|z_0|t -
   \cos|z_0|t\cot\frac{|z_0|t_0}{2} \right)[x_0,Jx_0] &&\\
{}-\frac{1}{2|z_0|^3}\sin|z_0|t\,\Big( \cos|z_0|t_0 - \cos|z_0|t -
   \cos|z_0|(t-t_0) + 1\Big)\,[J_\zeta x_0,Jx_0] &&\\
{}-\frac{1}{2|z_0|^3}\cos|z_0|t\,\Big( \sin|z_0|t - \sin|z_0|t_0 -
   \sin|z_0|(t-t_0)\Big)\,[JJ_\zeta x_0,x_0] &&\\
{}-\frac{1}{2|z_0|^2}\cos|z_0|t\,\Big( \cos|z_0|t_0 - \cos|z_0|t -
   \cos|z_0|(t-t_0) + 1\Big)\,[J_\zeta x_0,x_0] &&\\
{}-\frac{1}{2|z_0|^4}\sin|z_0|t\,\Big( \sin|z_0|t - \sin|z_0|t_0 -
   \sin|z_0|(t-t_0)\Big)\,[JJ_\zeta x_0,Jx_0] &&\\
{}=\, cz_0\, +\! &\zeta & \!.
\end{eqnarray*}

Using Lemma \ref{l3} and the identities we derived from it earlier, this 
becomes
\begin{eqnarray*}
\lefteqn{\dz(t)\, =\, }\\
&& c \left\{ 1 - \frac{\langle x_0,x_0\rangle t_0}{2|z_0|} \left(
   \cot\frac{|z_0|t_0}{2}-\sin|z_0|t-\cos|z_0|t\cot\frac{|z_0|t_0}{2}
   \right)\right\} z_0 \\
&& {}+\left\{ 1+\frac{\langle x_0,x_0\rangle}{2|z_0|^2}\,\Big( 1 -
   \cos|z_0|(t-t_0)+\cos|z_0|(2t-t_0)-\cos|z_0|t \Big)\right\}\zeta\\
&& {}+\frac{1}{2|z_0|^3}\,\Big\{ \sin|z_0|(t-t_0) - \sin|z_0|(2t-t_0) +
   \sin|z_0|t\Big\}\,[JJ_\zeta x_0,x_0]\,.
\end{eqnarray*}
Integrating both sides under $z(0)=0$ yields
\begin{eqnarray*}
z(t)
&=& c\,\bigg\{ t-\frac{\langle x_0,x_0\rangle t_0}{2|z_0|}\,\bigg(
   \cot\frac{|z_0|t_0}{2}t + \frac{1}{|z_0|}\cos|z_0|t \\
&& \hspace{5em}{}- \frac{1}{|z_0|} -
   \frac{1}{|z_0|}\sin|z_0|t\cot\frac{|z_0|t_0}{2} \bigg)\bigg\}\,z_0\\
&& {}+\bigg\{ t + \frac{\langle x_0,x_0\rangle}{2|z_0|^2}\,\bigg( t -
   \frac{1}{|z_0|}\sin|z_0|(t-t_0) + \frac{1}{2|z_0|}\sin|z_0|(2t-t_0) \\
&& \hspace{5em}{}-\frac{1}{|z_0|}\sin|z_0|t -\frac{1}{2|z_0|}\sin|z_0|t_0
   \bigg) \bigg\}\,\zeta\\
&& {}+\frac{1}{2|z_0|^3}\,\bigg\{ -\frac{1}{|z_0|}\cos|z_0|(t-t_0) +
\frac{1}{2|z_0|}\cos|z_0|(2t-t_0) \\
&& \hspace{5em}{}- \frac{1}{|z_0|}\cos|z_0|t +
   \frac{1}{2|z_0|}\cos|z_0|t_0 + \frac{1}{|z_0|}\bigg\}\,[JJ_\zeta
   x_0,x_0]\,.
\end{eqnarray*}
Since $z(t_0)=0$, we find
\begin{eqnarray*}
z(t_0) &=& \frac{ct_0}{|z_0|^2}\,\bigg\{ \langle z_0,z_0\rangle + \langle
   x_0,x_0\rangle - \frac{\langle x_0,x_0\rangle|z_0|t_0}{2}
   \cot\frac{|z_0|t_0}{2} \bigg\}\,z_0 \\
&& {}+\frac{1}{2|z_0|^2}\bigg\{ \Big(\langle\dg,\dg\rangle+\langle
z_0,z_0\rangle\Big)t_0 - \frac{\langle
x_0,x_0\rangle}{|z_0|}\sin|z_0|t_0 \bigg\}\,\zeta\, =\, 0\,.
\end{eqnarray*}

If $\zeta=0$ and $c\neq 0$, then
$$ \langle \dg,\dg\rangle = \langle x_0,x_0\rangle
\frac{|z_0|t_0}{2}\cot\frac{|z_0|t_0}{2} \,.$$
If $\zeta\neq 0$ and $c=0$, then
$$ |z_0|t_0 = \frac{\langle x_0,x_0\rangle}{\langle \dg,\dg\rangle+\langle
z_0,z_0\rangle}\sin|z_0|t_0\,.$$
If $\zeta\neq 0$ and $c\neq 0$, then
$$ \langle\dg,\dg\rangle = \langle
x_0,x_0\rangle\frac{|z_0|t_0}{2}\cot\frac{|z_0|t_0}{2}\,,$$
$$ |z_0|t_0 = \frac{\langle x_0,x_0\rangle}{\langle \dg,\dg\rangle+\langle
z_0,z_0\rangle}\sin|z_0|t_0\,.$$
Thus we have shown that if $\g(t_0)$ is a conjugate point along
$\gamma$ and $ t_0 \notin (2\pi/|z_0|)\bZ^*$, then $t_0\in\A_1 \cup \A_2$.

Note that $u$ and $z$ are uniquely determined by the constants $c$ and 
$\zeta$. Thus we have these multiplicities:
\begin{enumerate}
\item if $t_0 \in \A_1-\A_2$, then $\mcp(t_0) = 1$;
\item if $t_0 \in \A_2-\A_1$, then $\mcp(t_0) = \dim\z-1$;
\item if $t_0 \in \A_1\cap\A_2$, then $\mcp(t_0) = 1+\dim\z-1=\dim\z$.
\end{enumerate}

\medskip

The proof of the second part of Theorem \ref{th1} is similar.

\bigskip

For the third part, let $ Y(t) = z(t) + e^{tJ}u(t)$ be a Jacobi field
with $Y(0) = Y(t_0 )=0$.  Then by Proposition \ref{yj} we have
\begin{equation}
\dz(t) - [e^{tJ}u(t) , x'(t)] = cz_0 + d\zeta +\xi\, ,\label{e1}
\end{equation}
where $c,d$ are real constants and $\zeta ,\xi$ are constants in $\z$
such that $\langle z_0 , \zeta \rangle =1$, $ \langle \zeta , \zeta
\rangle = 0$, $\langle z_0 , \xi \rangle = \langle \zeta , \xi \rangle =
0$, and
\begin{equation}
\label{e2}
e^{tJ}\ddu(t) + e^{tJ} J \du(t) -J_{cz_0 + d\zeta +\xi }x'(t) = 0\,.
\end{equation}
Since $J^2 = -\langle z_0 , z_0 \rangle I = 0$, (\ref{e2}) reduces to
\begin{equation}
\label{te2}
(I + tJ ) \ddu + J\du = c J x_0 + dJ_\zeta x_0 + J_\xi  x_0
   +(td) J_\zeta J x_0 + tJ_\xi  Jx_0 \,.
\end{equation}
Note that the solution $u$ of (\ref{te2}) satisfying $u(0)=u(t_0)=0$ is
unique (if it exists).

Consider the vector field $u$ along $\g$ given by
$$ u(t) = \alpha_0 x_0 + \alpha_1 Jx_0 + \alpha_2 J_\zeta x_0 + \alpha_3
JJ_\zeta x_0 + \beta_1 J_\xi x_0 +\beta_2 JJ_\xi x_0 $$
where
\begin{eqnarray*}
\alpha_0 &=& {\txs\frac{d}{3}}t(t_0^2 - t^2 )\, ,\\
\alpha_1 &=& {\txs\frac{c}{2}}t(t-t_0) + {\txs\frac{d}{12}}t(3t^3 -2
   t^2_0t-t_0^3)\, ,\\
\alpha_2 &=& {\txs\frac{d}{2}}t(t-t_0)\, ,\\
\alpha_3 &=& {\txs\frac{d}{4}}t(t_0 +2t )(t_0 -t)\, ,\\
\beta_1 &=&\half t(t-t_0 )\, ,\\
\beta_2 &=& \quar t(t_0 +2t)(t_0 -t)\,.
\end{eqnarray*}
By a direct computation using Lemma \ref{l3}, this $u$ is a solution of
(\ref{te2}). Since $u(0)=u(t_0)=0$, it is the unique such solution.

Now we compute
\begin{eqnarray*}
[e^{tJ}u(t), e^{tJ}x_0 ] &=& [(I+tJ)u ,  (I+tJ)x_0 ]\\
&=& [u, x_0]+t[Ju,x_0] +t[u,Jx_0]+t^2[Ju,Jx_0]\\
&=& \alpha_1 [Jx_0 , x_0] + \alpha_2 [J_\zeta x_0 , x_0]+\alpha_3
   [JJ_\zeta x_0 ,x_0]\\
&& {}+ \beta_1 [J_\xi x_0, x_0] +\beta_2 [ JJ_\xi x_0, x_0]\\
&& {}+t\Big( \alpha_0 [Jx_0 ,x_0 ]+\alpha_2[JJ_\zeta x_0 ,x_0]+\beta_1
   [JJ_\xi x_0 , x_0]\Big)\\
&& {}+t\Big( \alpha_0 [x_0 ,Jx_0 ]+\alpha_2[J_\zeta x_0 ,Jx_0]+\alpha_3
   [JJ_\zeta x_0 ,J x_0]\\
&& {}+\beta_1 [J_\xi x_0 , Jx_0]+\beta_2 [JJ_\xi x_0 ,Jx_0]\Big)\\
&& {}+t^2 \Big(\alpha_2 [JJ_\zeta x_0, Jx_0 ] +\beta_1 [JJ_\xi x_0,
   Jx_0]\Big)\,.
\end{eqnarray*}
Observe that
\begin{eqnarray*}
[JJ_\zeta x_0,x_0] &=& [J_\zeta x_0,Jx_0]\,,\\
{}[JJ_\xi x_0,x_0] &=& [J_\xi x_0,Jx_0]\,,\\
{}[JJ_\zeta x_0,Jx_0] &=& [(-J_\zeta Jx_0-2I)x_0,Jx_0]\\
   &=& \langle Jx_0,Jx_0\rangle\zeta-2[x_0,Jx_0]\\
   &=& -2\langle x_0,x_0\rangle z_0\,,
\end{eqnarray*}
and
\begin{eqnarray*}
[JJ_\xi x_0,Jx_0] &=& -[J_\xi Jx_0,Jx_0]\\
   &=& [Jx_0,J_\xi Jx_0]\\
   &=& \langle Jx_0,Jx_0\rangle\xi\\
   &=& 0\,.
\end{eqnarray*}
Thus
\begin{eqnarray*}
\lefteqn{[e^{tJ}u(t), e^{tJ}x_0 ]}\qquad\quad\\
&=& -\langle x_0,x_0\rangle\alpha_1z_0 - \langle
   x_0,x_0\rangle\alpha_2\zeta + (\alpha_3+2t\alpha_2)[JJ_\zeta x_0,x_0]\\
&& {}- \langle x_0,x_0\rangle\beta_1\xi + (\beta_2+2t\beta_1)[JJ_\xi
   x_0,x_0]\\
&& {}- 2t\alpha_3\langle x_0,x_0\rangle z_0 - 2t^2\alpha_2 \langle
   x_0,x_0\rangle z_0 \\
&=& -\langle x_0,x_0\rangle(\alpha_1 + 2t\alpha_3 + 2t^2\alpha_2)z_0
   - \langle x_0,x_0\rangle\alpha_2\zeta - \langle
   x_0,x_0\rangle\beta_1\xi\\
&& {}+ (\alpha_3 + 2t\alpha_2)[JJ_\zeta x_0,x_0]
+ (\beta_2 + 2t\beta_1)[JJ_\xi x_0,x_0]\,.
\end{eqnarray*}
{F}rom this and (\ref{e1}), we get
\begin{eqnarray*}
\dz(t) + \langle x_0,x_0\rangle(\alpha_1 + 2t\alpha_3+2t^2\alpha_2)z_0 
   +\langle x_0,x_0\rangle\alpha_2\,\zeta && \\
{}+ \langle x_0,x_0\rangle\beta_1\,\xi 
   + (\alpha_3+2t\alpha_2)[JJ_\zeta x_0,x_0] && \\
{}+ (\beta_2+2t\beta_1)[JJ_\xi x_0,x_0] &=& cz_0+d\zeta+\xi\,.
\end{eqnarray*}
Integrating under $z(0)=0$,
\begin{eqnarray*}
z(t) + \langle x_0,x_0\rangle \bigg\{ \frac{c}{2}(\third t^3-\half t^2t_0) 
   + d\left(\frac{1}{20}t^5 - \frac{1}{8}t^4t_0 + \frac{1}{9}t^3t_0 - 
   \frac{1}{24}t^2t_0^3\right)\!\bigg\}\,z_0 && \\
{}+\frac{d\langle x_0,x_0\rangle }{2}\left(\third t^3 - \half 
   t^2t_0\right)\zeta + \frac{\langle x_0,x_0\rangle }{2}\left(\third t^3
   - \half t^2t_0\right)\xi \\
{}+\frac{d}{8}t^2(t-t_0)^2[JJ_\zeta x_0,x_0] + \frac{1}{8}t^2(t-t_0)^2 
   [JJ_\xi x_0,x_0] \\
=\, \left(cz_0+d\zeta+\xi\right)t\hspace{-1em} &&
\end{eqnarray*}
Using $z(t_0)=0$, from this we get
$$ -\langle x_0,x_0\rangle\left(\frac{c}{12}t_0^3 + 
\frac{d}{180}t_0^5\right)z_0 - \frac{d}{12}\langle x_0,x_0\rangle 
t_0^3\,\zeta - \frac{\langle x_0,x_0\rangle }{12}t_0^3\,\xi = 
\left(cz_0+d\zeta+\xi\right)t_0 $$
which implies that $d=0$ and $-\langle x_0,x_0\rangle t_0^2/12 = 1$. Thus
\begin{eqnarray*}
u(t) &=& \frac{c}{2}t(t-t_0)Jx_0 + \half t(t-t_0)J_\xi x_0 + \quar 
t(t_0+2t)(t_0-t)JJ_\xi x_0\,,\\
z(t) &=& ct\left\{ 1-\langle x_0,x_0\rangle t(\sixth t-\quar 
   t_0)\right\}z_0 + t\left\{ 1-\langle x_0,x_0\rangle t(\sixth t-\quar 
   t_0)\right\}\xi \\
&&{}- \frac{1}{8}t^2(t-t_0)^2[JJ_\xi x_0,x_0]\,.
\end{eqnarray*}
Note that $u$ and $z$ are uniquely determined by $c$ and $\xi$. Hence 
$\mcp(t_0) = \dim\z-1$.

\section{Proof of Theorem \protect\ref{th2}}
\label{pf2}

\begin{lemma}
Let\label{l1} $N$ be a pseudoriemannian, 2-step nilpotent Lie group with
nondegenerate center.  If $[x,J_z x]=cz$ for some real $c$ for every
nonnull $x\in\v$ and $z\in\z$, then there exists a selfadjoint operator
$A$ on $\v$ such that $J_z^2=\langle z,z\rangle A$ for all $z\in\z$
\end{lemma}
\begin{proof}
{\em Step 1:} $[x,J_z x]=cz$ for every $x\in\v$ and every nonnull
$z\in\z$.  Suppose $|x|=0$.  Then there exists a null $v\in\v$ such that
$\langle x,v\rangle\neq 0$.  Note $x+av$ is nonnull for every nonzero
scalar $a$.  Thus \begin{eqnarray} [x+av,J_z(x+av)] &=& c_a z \label{j1}\\
{}[x-av,J_z(x-av)] &=& c_{-a}z \label{j2} \end{eqnarray} where $c_{\pm a}$
are constants.  From these equations, we get
$$ [x,J_z x] + a^2[v,J_z v] =\half(c_a + c_{-a})z\,. $$
Since $a$ is arbitrary,
$$ [x,J_z x]\in\lsp z\rsp\,. $$

\noindent{\em Step 2:} $J_z J_{z'}+J_{z'}J_z=0$ for every mutually
orthogonal pair of nonnull vectors $z,z'$ in $\z$. By Step 1, we have
$$ [x+v,J_z(x+v)]=cz $$
for all $x,v\in\v$. Thus
$$ \langle J_{z'}(x+v),J_z(x+v)]=0 \,. $$
This implies
$$ \langle J_{z'}x,J_z x\rangle + \langle J_{z'}x,J_z v\rangle + \langle
J_{z'}v,J_z x\rangle + \langle J_{z'}v,J_z v\rangle = 0\, ,$$
so
$$ \langle -J_zJ_{z'}x,x\rangle + \langle -J_zJ_{z'}x,v\rangle + \langle
   -J_{z'}J_z x,v\rangle + \langle -J_zJ_{z'}v,v\rangle = 0\,. $$
Now
\begin{eqnarray*}
\langle -J_zJ_{z'}x,x\rangle &=& \langle -z,[J_{z'}x,x]\rangle\\
   &=& \langle -z,dz'\rangle\\
   &=& 0
\end{eqnarray*}
for some scalar $d$. Substituting this result into the previous equation 
yields
$$ \langle (J_zJ_{z'}+J_{z'}J_z)x,v\rangle = 0\, ,$$
which implies the desired conclusion.

\smallskip\noindent{\em Step 3:} $J_{z}^2=J_{z'}^2$ if $0\neq \langle
z,z\rangle = \langle z',z'\rangle$ and $z\perp z'$. By Step 1,
$$ [x,J_{z+z'}x]=c(z+z') $$
for all $x\in\v$. This implies
$$ \langle J_{z-z'}x,J_{z+z'}x\rangle=0 $$
whence
$$ \langle J_{z}x,J_z x\rangle - \langle J_{z'}x,J_{z}x\rangle + \langle
   J_z x,J_{z'}x\rangle - \langle J_{z'}x,J_{z'}x\rangle = 0 $$
so
$$ \langle J_z^2 x,x\rangle = \langle J_{z'}^2 x,x\rangle\,. $$
Thus, for all $x,v\in\v$,
$$ \langle J_z^2(x+v),x+v\rangle = \langle J_{z'}^2(x+v),x+v\rangle $$
whence
$$ \langle J_z^2 x,v\rangle = \langle J_{z'}^2 x,v\rangle $$
and therefore
$$ J_z^2 = J_{z'}^2\,. $$

\noindent{\em Step 4:} $J_{z}^2=-J_{z'}^2$ if $0\neq \langle
z,z\rangle = -\langle z',z'\rangle$ and $z\perp z'$. By Step 1, for all
$x\in\v$,
\begin{eqnarray*}
[x,J_{z+\frac{1}{2}z'}x] &=& c_1(z+\half z')\, ,\\
{}[x,J_{z}x] &=& c_2 z\, ,\\
{}[x,J_{\frac{1}{2}z'}x] &=& c_3\half z'\,.
\end{eqnarray*}
{}From these we infer
\begin{eqnarray*}
\langle J_z x,J_z x\rangle &=& c\langle z,z\rangle \\
\langle J_{z'}x,J_{z'}x\rangle &=& c\langle z',z'\rangle
\end{eqnarray*}
so
\begin{eqnarray*}
\langle -J_{z}^2 x,x\rangle + \langle -J_{z'}^2 x,x\rangle &=& 0\\
\langle \left(J^2_z + J^2_{z'}\right)x,x\rangle &=& 0
\end{eqnarray*}
and therefore
$$ J_z^2 = -J_{z'}^2\,. $$

\noindent{\em Step 5:} for an orthonormal basis $\{z_i\}$ of $\z$ such
that $\langle z_1,z_1\rangle = \cdots = \langle z_l,z_l \rangle = 1$ and
$\langle z_{l+1},z_{l+1}\rangle = \cdots = \langle z_k,z_k\rangle = -1$
and $J^2_{z_1} = A$, then $J^2_z = \langle z,z\rangle A$ for all $z\in\z$.
Expand $z$ in the basis $\{z_i\}$ as $z = \sum_i a_i z_i$. Then
\begin{eqnarray*}
J^2_z &=& \sum_i a_i J_{z_i}\,\sum_i a_i J_{z_i}\\
&=& \sum_i a_i^2 J_{z_i}^2 \qquad\mbox{by Step 2}\\
&=& \sum_{i=1}^l a_i^2 J_{z_1}^2 - \sum_{i=l+1}^k a_i^2 J_{z_1}^2
   \qquad\mbox{by Steps 3 and 4}\\
&=& \left( \sum_{i=1}^l a_i^2 - \sum_{i=l+1}^k a_i^2 \right) J_{z_1}^2\\
&=& \langle z,z\rangle A\,.\eop
\end{eqnarray*}
\end{proof}

\begin{lemma}
Let\label{l2} $N$ be pseudoregular.  For a nonzero $x\in\v$ and a nonnull
$z\in\z$, if $x,J_z x, J^2_z x$ are linearly dependent and $x,J_z^2 x$ are
linearly independent, then $x$ is null.
\end{lemma}
\begin{proof}
By assumption, we can find constants $a,b,c$, not all zero, such that
$ax+bJ_z x+cJ_z^2 x = 0$, and we may assume that $a\neq 0\neq b$ since $x$
and $J^2_z x$ are linearly independent. Since
$$ \ker\left( aI + bJ_z + cJ_z^2\right) \cap \ker\left( aI - bJ_z +
   cJ_z^2\right) = 0\, ,$$
the map
$$ aI-bJ_z+cJ_z^2 : \ker\left( aI+bJ_z+cJ_z^2\right)\to \ker\left(
   aI+bJ_z+cJ_z^2\right) $$
is nonsingular. Thus we may choose $v\in\ker\left(aI + bJ_z +
cJ_z^2\right)$ such that $ \left(aI-bJ_z+cJ_z^2\right)v=x$, so
\begin{eqnarray*}
\langle x,x\rangle &=& \langle\left(aI-bJ_z+cJ_z^2\right)v,x\rangle\\
&=& \langle v,\left(aI+bJ_z+cJ_z^2\right)x\rangle\\
&=& \langle v,0\rangle\\
&=& 0\,.\eop
\end{eqnarray*}
\end{proof}
\begin{proof}[Theorem \ref{th2}]
We first prove the ``only if" half, proceeding by steps and cases.
{\em Step 1:} $[x,J_z x]=cz$ for all nonnull $z\in\z$, $x\in\v$.
Suppose the contrary. Then there exist nonnull $z\in\z$, $x\in\v$ such
that $[x,J_z x]\notin\lsp z\rsp$. Consider the geodesic $\g$ such that
$\g(0)=1\in N$ and $\dg(0)=z+x$. Then there exists a connected, totally
geodesic, 3-dimensional submanifold $H$ with $\dg(0)\in\h=T_1H$ and
$\h\cap\z$ nonnull. Choose $z'\in\h\cap\z$ with $\langle z',z'\rangle\neq
0$. Then by the Gauss equation, we have $R(x+z,z')z' = -\quar J_{z'}^2 x
\in\h$ or $J_{z'}^2 x\in\h$. Similarly, $J_{z'}^4 x\in\h$.

{\em Case 1.1:} $J^2_{z'}x$ and $J^4_{z'}x$ are linearly independent. Then
it follows that $\h=\lsp z',J^2_{z'}x,J^4_{z'}x\rsp$.  This implies that
$x,z\in\h$ and $z$ and $z'$ are linearly dependent.  By the Gauss
equation, $R(x,z)x = -\quar[x,J_z x]\in\h$.  This and $\h=\lsp
z,J^2_{z}x,J^4_{z}x\rsp$ imply that $[x,J_{z}x]\in\lsp z\rsp$, a
contradiction.

\smallskip
{\em Case 1.2:} $J^2_{z'}x$ and $J^4_{z'}x$ are linearly dependent. This
and the nonsingularity of $J^2_{z'}$ imply that $x$ and $J^2_{z'}x$ are
linearly dependent, whence $x,z\in\h$.

{\em Subcase 1.2.1:} $z$ and $z'$ are linearly independent. Then
$\h = \lsp x,z,z'\rsp$. By the Gauss equation, $R(x,z)z =
-\quar J^2_z\in\h$ so $J^2_z x=cx$ for some constant $c$, and $R(x,z)z' =
-\quar J_zJ_{z'}x\in\h$ so $J_zJ_{z'}x = dx$ for some constant $d\neq 0$.
Thus $\add{x}{(z-\frac{c}{d}z')} = 0$.  This implies $\langle [x,v],
z-\frac{c}{d} z'\rangle = 0$ for all $v\in\v$.  But this contradicts $N$
being pseudoregular.

{\em Subcase 1.2.2:} $z$ and $z'$ are linearly dependent. Then $\h=\lsp
x,z,z_1+x_1\rsp$ for some $z_1+x_1\in\z\ds\v$. If $x$ and $x_1$ are
linearly dependent, then $\h=\lsp x,z,z_1\rsp$ and an argument similar to
that of the preceding subcase yields a contradiction. Thus we assume that
$x$ and $x_1$ are linearly independent. By the Gauss equation, $R(x,z)x =
-\quar [x,J_zx]\in\h$ which implies that $[x,J_zx]\in\lsp z\rsp$.

Therefore $[x,J_zx]\in\lsp z\rsp$ for all nonnull $z\in\z$, $x\in\v$.

\medskip
Step 1 and Lemma \ref{l1} imply that $J_z^2 = \langle z,z\rangle A$ for
all $z\in\z$ for some fixed, selfadjoint operator $A$ on $\v$.

\smallskip\noindent{\em Step 2:} $A=\l I$ for a nonzero constant $\l$.
Suppose the contrary; then we can choose nonnull $z\in\z$ and $x\in\v$
such that $x$ and $J_z^2 x$ are linearly independent.

Consider a geodesic $\g$ emanating from $1\in N$ with $\dg(0)=z+x$.  Then
there exists a totally geodesic submanifold $H$ with $\dg(0)\in T_1H=\h$
and $\h\cap\z$ nonnull.  Choose a nonnull $z'\in \h\cap\z$.  By the Gauss
equation, $J^2_{z'}x,J^4_{z'}x\in\h$.  By Step 1, $J^2_{z'}x = cJ^2_{z}x$
and $J^4_{z'}x = c^2 J^4_{z}x$.  Since $x$ and $J^2_{z}x$ are linearly
independent, so are $J^2_{z'}x$ and $J^4_{z'}x$.  Also, $\dg(0) = z+x \in
\h$ whence $\h = \lsp z, J^2_zx, J^4_zx\rsp = \lsp z,x,J^2_zx\rsp$.  Using
$J^2_z = \langle z,z\rangle A$, it follows that $[x,J^2_zx]=0$ for all
$z\in\z$.

Now consider a geodesic $\g$ emanating from $1\in N$ with $\dg(0)=x$. Then
$\g$ lies in $H$, $\g(t) = \exp(tx)$, $x$ and $J^2_zx$ are parallel along
$\g$, and so there exists a parallel vector field $\alpha z + \beta J_zx$
along $\g$ where $\alpha(0) = 1$ and $\beta(0)=0$.

{\em Case 2.1:} $[x,J_z x]=cz$ for $c=l^2$ with $l>0$. Then
$(\cos\frac{lt}{2})z + (\frac{1}{l}\sin\frac{lt}{2})J_zx$ is the desired
parallel vector field. Applying the Gauss equation at the point $\g(t) =
\exp(tx)$, we obtain
\begin{eqnarray*}
R(x,\alpha z+\beta J_zx)x &=& \alpha R(x,z)x + \beta R(x,J_zx)x\\
&=& -\quar\alpha[x,J_zx] + \tquar\beta J_{[x,J_zx]}x\\
&=& -\quar c\alpha z + \tquar c\beta J_zx \in\lsp x,\alpha z+\beta
J_zx,J^2_zx\rsp\,.
\end{eqnarray*}
Thus we can write $-\quar c\alpha z + \tquar c\beta J_zx = s_1x+s_2(\alpha
z+\beta J_zx)+s_3 J^2_zx$ for suitable constants $s_i$. Then $s_2=-\quar
c$ and $s_1x-c\beta J_zx+s_3 J^2_zx=0$. Since $x$ and $J^2_zx$ are
linearly independent, Lemma \ref{l2} implies that $x$ is null, which is a
contradiction. Thus this case cannot happen.

{\em Case 2.2:} $[x,J_z x]=cz$ for $c=-l^2$ with $l>0$. Then
$(\cosh\frac{lt}{2})z + (\frac{1}{l}\sinh\frac{lt}{2})J_zx$ is parallel
along $\g$. By a similar argument to that of the previous case, we have
another contradiction.

{\em Case 2.3:} $[x,J_z x]=0$. This condition implies $\langle
x,J^2_zx\rangle=0$.  Since $x$, $J^2_zx$, and $J^4_zx$ are linearly
dependent, there exist constants $0\neq a,b,c\neq 0$ such that
$ax+bJ^2_zx+cJ^4_zx=0$.  This implies that $\langle ax + bJ^2_zx +
cJ^4_zx,x\rangle=0$ whence $a\langle x,x\rangle +c\langle
J^4_zx,x\rangle=0$.  Since $\langle x,x\rangle\neq0$, $\langle
J^4_zx,x\rangle\neq0$.

Now $z+\half tJ_zx$ is a parallel vector field along $\g$.  Applying the
Gauss equation,
\begin{eqnarray*}
\lefteqn{R(J^2_zx,z+\half tJ_zx)J^2_zx }\qquad\ \:\quad\\
&=& -\quar[J^2_zx,J_zx] + {\txs \frac{3}{8}}t\add{x}{[J^2_zx,J_zx]}
\in \lsp x,J_z^2 x, z+\half tJ_z x\rsp\, ,
\end{eqnarray*}
and we know that $R(J^2_zx,z)J^2_zx = -\quar[J^2_zx,J_zx]\in\h$.  This
last implies that $[J^2_zx,J_zx]=dz$ for some constant $d$.  Since
$\langle x,J^4_zx\rangle\neq 0$, we have $d\neq 0$.  Thus
$$-\quar dz + {\txs \frac{3}{8}}tdJ_zx\in\lsp x,J^2_zx,z+\half
   tJ_zx\rsp$$
so there exist constants $s_i$ such that
$$-\quar dz+{\txs \frac{3}{8}}tdJ_zx = s_1x+s_2J^2_zx+s_3(z+\half
   tJ_zx)\,,$$
whence $s_3=-\quar d$ and $s_1x-\half tdJ_zx+s_2J^2_zx=0$. But the last
equation contradicts Lemma \ref{l2}.

Therefore we conclude that $A=\l I$. The ``only if'' half now follows, 
$-\l$ being the required constant of homothety.

\medskip
Conversely, assume that $N$ is homothetic to one of pseudo$H$-type and
consider a geodesic $\g$ with $\g(0)=1\in N$ and $\dg(0)=z_0+x_0$.  We may
as well assume $|z_0|\neq 0$.

If $x_0=0$, then choose a nonzero $x\in\v$ such that $x$ and $Jx$ are
linearly independent. It follows that $H=\exp\lsp z_0,x,Jx\rsp$ is a
3-dimensional, totally geodesic subgroup of $N$ with $\dg(0)\in T_1H$ as
desired.

Now assume $x_0\neq 0$. If $x_0$ and $Jx_0$ are linearly independent then,
almost as before, $H = \exp\lsp z_0, x_0, Jx_0\rsp$ is a suitable
subgroup. If $x_0$ and $Jx_0$ are linearly dependent, then there exists a
nonzero scalar $d$ such that $Jx_0=dx_0$. Now it follows that $J^2 =
d^2I$, so $\v$ can be decomposed as $\v = \ker(J-dI)\ds\ker(J+dI)$.  One
may show (as in \cite{JP2}, for example) that $\ker(J-dI)$ and
$\ker(J+dI)$ are complementary null subspaces of $\v$.  Thus we may choose
$y\in\ker(J+dI)$ such that $\langle x_0,y\rangle = 1$.  Letting $\{z_0,
z_1, \ldots , z_k\}$ be an orthogonal basis of $\z$ and using
$JJ_{z_i}+J_{z_i}J=0$ for $i\neq 0$ from Lemma \ref{l3}, we find $\{ J_z
x_0 \mid z\in\lsp z_1, \ldots , z_k\rsp\} \subseteq \ker(J+dI)$.  So
$\langle J_{z_i}x_0,y\rangle = 0$ for $1\le i\le k$ since
$J_{z_i}x_0,y\in\ker(J+dI)$.  Therefore $[x_0,y] = (d/\langle
z_0,z_0\rangle )z_0$.  This implies that $H = \exp\lsp z_0,x_0,y\rsp$ is a
3-dimensional, totally geodesic subgroup with $\dg(0) \in T_1H$ as
required.
\end{proof}

\end{document}